\documentclass[12pt,a4paper]{article}
\usepackage{amsmath, amssymb, hyperref}
\usepackage{tikz}
\usepackage{tikz-cd}
\def\r){\right)}
\def\l({\left(}
\def\qmbox#1{\quad\mbox{#1}\quad}
\def\pmx#1{\begin{pmatrix}#1\end{pmatrix}}

\def\dsa{d^*\! A}

\def\R{\mathbb{R}}
\def\C{\mathbb{C}}
\def\H{\mathbb{H}}
\def\O{\mathbb{O}}

\def\RRR{\mathcal{R}}
\def\RR{\mathcal{R}^*}

\def\CP{\C\mathrm{P}}
\def\CH{\C\mathrm{H}}
\def\HP{\H\mathrm{P}}
\def\HH{\H\mathrm{H}}
\def\OP{\O\mathrm{P}}
\def\OH{\O\mathrm{H}}

\def\u{\mathfrak{u}}
\def\h{\mathfrak{h}}
\def\k{\mathfrak{k}}
\def\g{\mathfrak{g}}
\def\m{\mathfrak{m}}

\def\so{\mathfrak{so}}

\def\sp{\mathfrak{sp}}

\def\norm#1{\left\|#1\right\|}
\def\vol{\mathbf{vol}\,}

\DeclareMathOperator{\Hom}{Hom}
\DeclareMathOperator{\id}{id}

\DeclareMathOperator{\trace}{trace}

\def\st#1{\left\{#1\right\}}
\def\set#1#2{\left\{#1\ \vline\hspace{110sp}\vline\hspace{110sp}\vline\ #2\right\}}
\def\scp#1#2{\left\langle#1\ \vline\hspace{110sp}\vline\hspace{110sp}\vline\ #2\right\rangle}

\def\spn#1{\mathop{\left\langle#1\right\rangle}}
\def\omt#1{\widehat{#1}}

\def\cases#1{\left\{\begin{array}{c@{\text{ if }}l}#1\end{array}\right.}

\def\[#1{\begin{equation}\label{#1}}
\def\]{\end{equation}}

\def\proof{{\bf Proof: }}
\def\proofend{\hfill$\bullet$\par}

\newtheorem{thm}[equation]{Theorem}

\newtheorem{lm}[equation]{Lemma}

\newtheorem{ntt}[equation]{Notation}

\def\longhookrightarrow{\lhook\joinrel\relbar\joinrel\relbar\joinrel\rightarrow}
\def\longepiarrow{\relbar\joinrel\relbar\joinrel\twoheadrightarrow}

\def\DD#1{\left.\frac{D}{d#1}\right|_{#1=0}}

\def\dd#1{\left.\frac{d}{d#1}\right|_{#1=0}}

\def\ad{\mathrm{ad}}

\numberwithin{equation}{section}

\begin{document}\sloppy

\def\email#1{}
\def\affiliation#1{}

\title{Green's functions, Biot-Savart Operators and Linking Numbers on Negatively Curved Symmetric Spaces}

\author{Stefan Bechtluft-Sachs\\
  \texttt{stefan.bechtluft-sachs@mu.ie}
  \and
  Evangelia Samiou\\
  \texttt{samiou@ucy.ac.cy}
}

\maketitle
\begin{abstract}
We construct radial fundamental solutions for the differential form Laplacian on negatively curved symmetric spaces.
At least one of these Green's functions also yields a Biot-Savart Opearator, i.e. a right inverse of the exterior differential on closed forms with image in the kernel of the codifferential. Any Biot-Savart operator gives rise to a Gauss linking integral.
\end{abstract}

\begin{center}\it Keywords: Laplace-de Rham Operator, Biot-Savart operator, Rank one symmetric space, Radial kernel function, Levinson Theorem, Linking number, Gauss linking integral\end{center}

\section{Introduction}
The purpose of this paper is to extend the well known radial analysis of the scalar Laplace operator on the Euclidean space $\R^m$ to the differential form Laplacian on negatively curved symmetric spaces.
In more detail, in $\R^m$, by rotational symmetry, Poisson's equation reduces to an ordinary differential equation.
For this reduction, the key notion is that of a radial kernel function.
For example, the electrostatic potential $\phi$ caused by a charge distribution $\rho$ must satisfy the Poisson equation
\[{coul}\Delta \phi=\rho\ .\]
If $\rho$ has compact support, then there is a radial solution
$$\phi(x)=[Q\rho](x)=\int_{\R^3}\frac{\rho(y)}{4\pi\norm{x-y}}\ d^3y\ .$$
Thus the scalar Laplacian $\Delta$ has a right inverse $Q$ with integral kernel $q(x,y)=A(d(x,y))$ depending on the distance only, and Poisson's equation becomes equivalent to the ordinary differential equation
\[{couloder3}\begin{split}
A''(t)+\frac{1}{t}A'(t)&=0\qmbox{with}\\
A(t)&\sim \frac{-1}{t\vol S^2}\text{ as }t\to 0\ ,
\end{split}\]
which has the solution $A(t)=\frac{1}{4\pi t}$.
This generalizes to harmonic spaces, i.e. Riemannian manifolds whose distance spheres have constant mean curvature.
All negatively curved symmetric spaces are harmonic (see section \ref{harmonicprelim}).

We will extend this construction of a fundamental solution to the differential form Laplace operator to negatively curved symmetric spaces $X$.
For a given compactly supported $l$-form, $\omega\in\Omega^l_c(X)$, we look for a solution $\mu\in\Omega^l(X)$ of $\Delta\mu=\omega$ in the form 
\[{intkerform}\mu_x=[Q\omega]_x=\int_X q_{x,y}\omega_y\ dy\ .\]
The difficulty here is that the kernel function $q$ of the integral operator $Q$ in \eqref{intkerform} consists of maps $q_{x,y}\colon\Lambda^lT_y^*X\to\Lambda^lT_x^*X$ relating different fibres of the bundle $\Lambda^lT^*X$ of alternating forms, i.e. in general different vector spaces.
To reduce the partial differential equation for the integral kernel $q$ in \eqref{intkerform} to an ordinary differential equation and fully exploit the two point symmetry of $X$, we need a reference map.
In this paper we will use the parallel transport, i.e. we will compare $q_{x,y}$ with the parallel transport along the unique geodesic from $y$ to $x$.
Our main Theorem \ref{thm:main} contains the radial equation \eqref{maindeq}, in ``parallel transport form'', in the terminology of \cite{dtg} and \cite{tk}, for the differential form Laplacian.
Radial Green's functions for the Laplacian on differential forms exist and are determined by solutions of a matrix valued ordinary differential equation resembling \eqref{couloder3} with the factor $\frac{1}{t}$ replaced by the mean curvature of the distance spheres and a more involved right hand side reflecting the more complicated curvature of $X$.

Hodge Theory provides a close link between Green's functions for the differential form Laplacian, Biot-Savart Operators and Gauss linking integrals on compact Riemannian manifolds.
In the non compact case this no longer holds.
Nonetheless, in Theorem \ref{thm:bs} we show that there is choice of a solution $A$ of the radial equation \eqref{maindeq} so that the right inverse $Q^A$ of $\Delta$ corresponding to $A$ gives a Biot-Savart operator with radial kernel.
In Theorem \ref{thm:gs} this leads to a Gauss integral formula for the linking number in negatively curved symmetric spaces.

For $3$-dimensional space forms, the papers \cite{ctg}, \cite{dtg} and \cite{tk} obtain explicit formulae for the differential form Green's function, the Biot-Savart Operator and Gauss linking integrals, by more direct means. The expressions of $q$ in terms of elementary functions are quite intricate, already in this case. In section \ref{hypsp} we derive the explicit radial equation for the Green's function in hyperbolic space and compare it with the result of \cite{dtg}.

Gauss integrals for the linking number do not depend on the full Biot-Savart operator, they do not require to solve the codifferential.
In space forms the linking number can also be expressed as a mapping degree.
The resulting explicit formulas are given in \cite{dtg3} and also \cite{kp}.

We are grateful to Marcos Salvai for initiating this work, for many profitable discussions during his visit to the University of Cyprus and for bringing to our attention the papers \cite{ctg}, \cite{dtg}, \cite{tk} and \cite{dtg3}.

\section{Preliminaries}\label{sec:prelim}
We recall some fundamental facts for the Laplace operator on differential forms and its inversion in $\R^m$ by means of the Newton potentials.
Then we discuss the consequences of fixing an axis in a rank one symmetric space, in particular the reduction of equivariant kernel functions to functions on $\R^+$. 

For the Laplace operator on forms we refer to \cite{B}, \cite{LM}. For symmetric spaces see \cite{Esch}, \cite{Helg}, \cite{Loos}. 
\subsection{The Weitzenb\"ock Formula}
The Laplace operator on differential forms on a Riemannian manifold $X^m$ is defined by
$$\Delta\omega = (dd^* + d^*\!d)\omega$$
where $d\colon\Omega^k(X)\to\Omega^{k+1}(X)$ denotes the exterior differential and $d^*$ its adjoint.
In terms of an orthonormal basis $e_1,\ldots,e_m$ of $T_pX$ these are given by
\[{diffcodiff}d\omega=e_\mu\wedge\nabla_{e_\mu}\omega\qmbox{,}d^*\!\omega=-e_\mu\neg\nabla_{e_\mu}\omega\]
where $\nabla$ is the Levi-Civita connection and for $\alpha\in\Omega^l(X)$,
$$e_\mu\wedge\alpha(x_0,x_1,\ldots,x_l)=\sum_{j=0}^l(-1)^j\scp{e_\mu}{x_j}\alpha(x_0,x_1,\ldots,\omt{x_j},\ldots,x_l)\ ,$$
$$e_\mu\neg\alpha(x_1,\ldots,x_{l-1})=\alpha(e_\mu,x_1,\ldots,x_{l-1})\ .$$
By the Weitzenb\"ock Formula
\[{WBF}\Delta=\nabla^*\nabla+\RR\]
where
$$\nabla^*\nabla\omega=-\trace\l(\nabla^2\omega\r)$$
and
$$\RR=-\sum_{\mu,\nu=1}^m e_\mu\wedge e_\nu\neg R_{e_\mu,e_\nu}\ .$$
Thus if $\omega\in\Omega^{l+1}(X)$, then $\RR\omega=\omega\circ\RRR$ with
\def\Ric{\mathrm{Ric}\,}
\begin{align*}
\RRR(v_0\wedge\ldots\wedge v_l)&=\sum_{i=0}^lv_0\wedge\ldots\wedge\underset{i}{\Ric{v_i}}\wedge\ldots\wedge v_l\\
&\hspace{-15mm}+\sum_{0\leq j<i\leq l}\ \sum_{\nu=1}^m v_0\wedge\ldots\wedge\underset{j}{R_{v_i,v_j}e_\nu}\wedge\ldots\wedge\underset{i}{e_\nu}\ldots\wedge v_l
\end{align*}
In low degrees $0$, $1$, $2$ this becomes
\begin{align*}
&l+1=0:&\RRR&=0\\
&l+1=1:&\RRR(v)&=\sum_\nu R_{v,e_\nu}e_\nu=\Ric(v)\\
&l+1=2:&\RRR(v\wedge w)&=\Ric(v)\wedge w+v \wedge\Ric(w)\\
&&&\hspace{5mm}+\sum_\nu R_{w,v}e_\nu \wedge e_\nu
\end{align*}

\subsection{The Newton Potential for differential forms on $\R^{m}$.}
\def\cinf{C^\infty(\R^{m})}
On $\R^m$ the Laplace Operator on forms takes a simple form.
For $0\leq l\leq m$, a $l$-form can be written as
$$\alpha=\sum_{1\leq i_1\leq\cdots\leq i_l\leq m} a_{i_1,\cdots i_l} dx_{i_1}\wedge dx_{i_2}\wedge\ldots\wedge dx_{i_l}$$
with $a_{i_1,\cdots i_l}\in\cinf$, where the $x_i$ denote the standard coordinates.
A straightforward calculation shows that the Laplace operator acts component-wise,
$$\Delta\alpha=\sum_{1\leq i_1\leq\cdots\leq i_l\leq m} (\Delta a_{i_1,\cdots i_l}) dx_{i_1}\wedge dx_{i_2}\wedge\ldots\wedge dx_{i_l}\ .$$
The Newton potentials are a Green's function for the Laplacian on functions.
They are defined by
\[{eq:newton}N(x,y):=\cases{\displaystyle\frac{1}{(m-2)\norm{x-y}^{m-2}\vol S^{m-1}} & m\neq 2 \\ \displaystyle\frac{-\log(\norm{x-y})}{2\pi} & m=2}\ .\]
If $\omega\in\Omega^l(\R^{m})$ is a differential $l$-form on $\R^{m}$ of compact support, then
$$\Delta_x\int_{\R^{m}}N(x,y)\omega_y\ d^{m}y=\omega\ .$$
If $\omega$ is closed then a solution $\alpha$ to the Cartan equation $d\alpha=\omega$ on $\R^m$ is given by,
\[{biotsavartrm}\alpha_x=d^*_x\int_{\R^{m}}N(x,y)\omega_y\ d^{m}y=\int_{\R^{m}}\frac{(x-y)\neg\omega_y}{\vol S^{m-1}\norm{x-y}^{m}}\ d^{m}y\ .\]
In the case $m=3$ this is the Biot-Savart formula for the magnetic field $\alpha$ caused by a stationary current $\omega$.

\subsection{Harmonic spaces}\label{harmonicprelim}
Harmonic spaces are Riemannian manifolds whose volume form in normal coordinates is a function of the distance from the origin. This is equivalent to all distance spheres of radius $t$ having constant mean curvature $h(t)$, a function of the radius only. A consequence is that the scalar Laplace operator has a Green's function $q(x,y)$ which is a radial kernel function in the sense of \cite{Szabo}, i.e. there is a function $A$ so that $q(x,y)=A(d(x,y)$.
Rank one symmetric spaces are harmonic. The converse, the famous Lichnerowicz conjecture \cite{Li}, was proved by Szabo \cite{Szabo} for compact spaces. There are however noncompact, non symmetric harmonic Riemannian manifolds (Damek Ricci spaces, \cite{DaRi}).

Let $h(t)$ be the mean curvature of the spheres of radius $t$ in a $m$-dimensional harmonic space $X$, $m>2$, and let $A\colon\R^+\to\R$ solve the ordinary differential equation
\[{coulode}\begin{split}
A''(t)+h(t)A'(t)&=0\qmbox{with}\\
A(t)&\sim \frac{1}{(m-2)t^{m-2}\vol S^{m-1}}\text{ as }t\to 0\ .
\end{split}\]
Then the integral operator $Q$ defined by
\[{coulsolform}[Q\rho](x)=\int_{X}A(d(x,y))\rho(y)\ dy\]
for compactly supported functions $\rho$ on $X$,
is a right inverse of the scalar Laplacian, i.e. $\Delta Q\rho=\rho$. In \eqref{coulsolform} and in the sequel integration in $X$ always is with respect to the Riemannian volume.
\medskip

\subsection{Negatively curved symmetric spaces}\label{rk1sym}
Negatively curved symmetric spaces are exactly the non compact rank one symmetric spaces. On the other hand rank one symmetric spaces are two point homogeneous, i.e. the isometry group acts transitively on equidistant pairs. In particular these spaces are harmonic. 

Let $X$ be a rank one symmetric space of dimension $m=n+1$.
We fix a point $p\in X$ and a unit tangent vector $T(p)\in ST_pX$.
Let $G$ be the group of isometries of $X$ and let $K=G_p$ be the isotropy group at $p$.
Let $H=G_{p,T}\subset K$ be the isotropy group of the vector $T$,
$$H=\set{g\in G}{gp=p\text{ and }dgT=T}\ .$$
We will identify the Lie algebra $\g$ of $G$ with that of the Killing fields of $X$.
The Lie algebras of $K$ and $H$ correspond to spaces of Killing fields,
\begin{align*}
\k&=\set{k\in\g}{k(p)=0}\qmbox{,}\\
\h&=\set{h\in\k}{\nabla_T k(p)=0}=\set{h\in\k}{[T,k]=0}
\end{align*}
where we have identified $T$ with its extension to a Killing field in
$$\m=\set{m\in\g}{\nabla m(p)=0}\ .$$
The scalar product on $\m$ will be so that the isomorphism
$$\m\to T_pX\qmbox{,}m\mapsto m(p)$$
is isometric. Let $\m_0$ be the orthogonal complement of $T$ in $\m$. Since $X$ is symmetric of rank one, the map
\[{eq:comp-h}\begin{array}{rcccl}
\k&\longhookrightarrow&\so(\m)&\longepiarrow&\m_0\\
k&\mapsto&\ad_k&\mapsto&[T,k]
\end{array}
\]
is surjective and has kernel $\h$.
If $X$ is not flat, i.e. not a euclidean space, then
$$\k_0:=[T,\m_0]$$
is a complement to $\h$ in $k$.
In case $X=\R^n$ we let $\k_0$ be the set of infinitesimal generators of rotations in planes containing $T$.
In any case we obtain splittings of $\g$, invariant under $Ad_K$ respectively $Ad_H$,
\[{liesplit}\g=\k\oplus\m=\underbrace{\h\oplus\k_0}_{=\k}\oplus\underbrace{\spn{T}\oplus\m_0}_{=\m\cong T_pX}\ , \]
and an isomorphism
$$\k_0\to\m_0\qmbox{,}k\mapsto[T,k]\ .$$
Below we list $\k$, $\h$ and the decomposition of $\m$ into irreducible $H$-modules for the rank one symmentric spaces.
\par
\def\spin{\mathfrak{spin}}
\begin{center}
\begin{tabular}{c|c|c|l}
$X$&$\k$&$\h$&$\m=\spn{T}\oplus\m_0$\\\hline
$\R^n$,$S^n$,$H^n$&$\so(n)$&$\so(n-1)$&$\spn{T}\oplus\R^{n-1}$\\\hline
$\CP^n$,$\CH^n$&$\u(n)$&$\u(n-1)$&$\spn{T}\oplus\spn{iT}\oplus\C^{n-1}$\\\hline
$\HP^n$,$\HH^n$&$\sp(1)\oplus\sp(n)$&$\sp(1)\oplus\sp(n-1)$&$\spn{T}\oplus\spn{iT,jT,kT}\oplus\H^{n-1}$\\\hline
$\OP^2$,$\OH^2$&$\so(9)=\spin(9)$&$\so(7)=\spin(7)$&$\spn{T}\oplus\R^7\oplus\O$\\
\end{tabular}
\end{center}
Because of the Jacobi identity, the operator $J\colon\m\to\m$, $J(m)=[T,[T,m]]$ is symmetric. In fact, for $n,m\in\m$, we have
\begin{align*}
\scp{[T,[T,m]]}{n}&=-\scp{T}{[n,[T,m]]}\\
&=-\scp{T}{[[n,T],m]}-\underbrace{\scp{T}{[T,[n,m]]}}_{=0}\\
&=-\scp{T}{[m,[T,n]]}
\end{align*}
because $[\m,\m]\subset\k$ and $\ad_k\in\so(\m)$ for all $k\in\k$.
\medskip

The geodesic $\gamma$ with $\gamma(0)=p$ and $\gamma'(0)=T$ coincides with the one parameter orbit generated by $T$, i.e. $\gamma(t)=e^{tT}p$ and $T$ is the generator of the one parameter group of translations along $\gamma$. 
At the point $p$, the operator $J$ coincides with the Jacobi operator along the geodesic, i.e. $Jm(p)=R_{T,m}T(p)$.
In the usual normalisation of the metric on $X$ it has eigenvalues $0$ if $X=\R^{m}$, $-1,-4$ if $X$ is of compact type and ${1,4}$ if $X$ is of non-compact type.

\begin{lm} Let $m\in\m$ be an eigenvector of $J$, $Jm=[T,[T,m]]=\lambda^2 m$ with $\lambda\neq 0$. Let $k=\frac{1}{\lambda^2}[T,m]\in\k$.
Then
\[{nablak}\nabla_Tk(p)=m(p)\ ,\]
\[{krot}[T,k]=m\qmbox{and} [m,k]=-T\ .\]
In particular, $\ad_k\colon x\mapsto [x,k]$ preserves the plane spanned by $T$ and $m$.
\end{lm}
\proof The first equation in \eqref{krot} holds by assumption, and since $[T,k]=\nabla_Tk-\nabla_kT$ and $k(p)=0$, this implies \eqref{nablak}.

It suffices to to show the second equation at the point $p$ since both sides are Killing fields in $\m$.
On the unit sphere $S\m\subset\m$ consider the function $S\m\to\R$, $x\mapsto\scp{[m,[x,m]]}{x}=\scp{R_{x,m}m}{x}$. At $x=T$ this takes an extreme value, hence $T$ is critical. It follows that $\scp{[m,[T,m]]}{x}=0$ whenever $x\perp T$, hence $[m,[T,m]]=\lambda^2 T$.
\proofend
\medskip

\begin{ntt}\label{ntt-basis} We fix an orthonormal basis $\set{m_i}{i=1\ldots n}$ of $\m_0$ consisting of eigenvectors of $J$, i.e. so that $Jm_i=\lambda_i^2 m_i$, (i.e. $R_{T,m_i}T=\lambda_i^2 m_i$).
If $X$ is not flat, we let $k_i=\frac{1}{\lambda_i^2}[T,m_i]\in\k_0$.
If $X=\R^m$, we let $k_i$ be the generator of the rotation in the plane spanned by $T$ and $m_i$.
In all cases we have a basis $\set{k_i}{i=1\ldots n}$ of $\k_0$ so that \eqref{nablak} and \eqref{krot} hold.
\end{ntt}

For $x,y\in X$ we denote by $P_{x,y}\colon T_yX\to T_xX$ the parallel transport from $y$ to $x$ along the minimizing geodesic.
Thus $P_{p,\gamma(t)}$ denotes the parallel transport along $\gamma$ from $\gamma(t)$ to $p=\gamma(0)$.
Since $X$ is symmetric, the Jacobi operator $V\mapsto R_{TV}T$ along $\gamma$ commutes with the parallel transport along $\gamma$, i.e.
$$R_{T,P_{\gamma(t),p}V}T=P_{\gamma(t),p}R_{TV}T=P_{\gamma(t),p}JV$$
for $V\in T_pX$.
The Killing fields $k\in\k$ are Jacobi fields along the geodesic $\gamma$ satisfying the second order differential equation
$$\nabla_T\nabla_Tk(\gamma(t))=R_{T,k(\gamma(t))}T=P_{\gamma(t),p}JP_{p,\gamma(t)}k\ .$$
If $m(p)=\nabla_Tk(p)$ is an eigenvector of $J$ corresponding to the eigenvalue $\lambda^2$, $\lambda\in i\R$ or $\lambda\in\R$, then
\[{eq:evaljac}
\begin{split}
k(\gamma(t))&=P_{\gamma(t),p}\frac{\sinh\l(\lambda t\r)}{\lambda}m(p)\qmbox{and}\\\nabla_Tk(\gamma(t))=m(\gamma(t))&=P_{\gamma(t),p}\cosh\l(\lambda t\r)m(p)\ .
\end{split}
\]
From these formulas we immediately compute the volume $\sigma(t)$ of the geodesic sphere of radius $t$ and its mean curvature $h(t)$,
\[{eq:sigmah}\begin{split}
\sigma(t)&
=\vol(S^n)\prod_{i=1}^n\frac{\sinh(\lambda_i t)}{\lambda_i}=t^n\l(\vol(S^n)+\frac{t^2}{6}\sum_{i=1}^n\lambda_i^2+\cdots\r)\\
h(t)&=\frac{\sigma'(t)}{\sigma(t)}=\sum_{i=1}^n\frac{\lambda_i\cosh(\lambda_i t)}{\sinh(\lambda_i t)}=\frac{n}{t}+\frac{t}{3}\sum_{i=1}^n\lambda_i^2+\cdots\ .
\end{split}\]
We will need to compute $\nabla_kk'$ for basis vectors $k,k'\in\k_0$ along $\gamma$.
\begin{lm} Let $m,m'\in\m_0$ be elements of the basis chosen in \ref{ntt-basis}, with corresponding $k,k'\in\k_0$ and eigenvalues $\lambda,\lambda'$. Then
$$\nabla_mk'(\gamma(t))=\cosh(\lambda t)\cosh(\lambda't)P_{\gamma(t),p}\ad_{k'}(m)(p)\ .$$
In particular,
\[{nablakinx}
\nabla_u k(\gamma(t))=\cosh(\lambda t)P_{\gamma(t),p}\ad_{k }P_{p,\gamma(t)}u
\]
for every $u\in T_{\gamma(t)}X$ and
\[{nablakk}\nabla_kk=-\frac{1}{\lambda}\sinh(\lambda t)\cosh(\lambda t) T\ .\]
\end{lm}
\proof Killing fields $X$ in a Riemannian manifold satisfy the differential equation
$$ R_{A,X} B=\nabla^2_{A,B}X$$
for any vector fields $A,B$. For $u,v,w\in\m$ we also have
$$R_{v,w}u(p)=[u,[v,w]](p)\ ,$$
see \cite{Esch} for example. By \eqref{eq:evaljac} and since $R$ is parallel, 
\begin{align*}
R_{T,k'}m(\gamma(t))&=R_{T,P_{\gamma(t),p}m'(p)\frac{\sinh(\lambda't)}{\lambda'}}P_{\gamma(t),p}m(p)\cosh(\lambda t)\\
&=\frac{\sinh(\lambda't)\cosh(\lambda t)}{\lambda'}P_{\gamma(t),p}R_{T,m'}m(p)\\
&=\frac{\sinh(\lambda't)\cosh(\lambda t)}{\lambda'}P_{\gamma(t),p}[m,[T,m']](p)\\
&=\lambda'\sinh(\lambda't)\cosh(\lambda t)P_{\gamma(t),p}[m,k']](p)\\ &=\lambda'\sinh(\lambda't)\cosh(\lambda t)P_{\gamma(t),p}\nabla_mk'(p)\ .
\end{align*}
On the other hand, by the differential equation for Killing fields, we can compute the same as
\begin{align*}
R_{T,k'}m(\gamma(t))&=\nabla^2_{T,m}k'(\gamma(t))\\
&=\nabla_T\nabla_mk'(\gamma(t))-\nabla_{\nabla_Tm}k'(\gamma(t))\\
&=\nabla_T\nabla_mk'(\gamma(t))-\frac{\lambda\sinh(\lambda t)}{\cosh(\lambda t)}\nabla_mk'(\gamma(t))\ .
\end{align*}
Thus $\nabla_mk'(\gamma(t))$ satisfies the differential equation
\begin{align*}
&\lambda'\sinh(\lambda't)\cosh(\lambda t)P_{\gamma(t),p}\nabla_mk'(p)\\
&\hspace{20mm}=\nabla_T\nabla_mk'(\gamma(t))-\frac{\lambda\sinh(\lambda t)}{\cosh(\lambda t)}\nabla_mk'(\gamma(t))\ .
\end{align*}
which has the unique solution
$$\nabla_mk'(\gamma(t))=\cosh(\lambda t)\cosh(\lambda't)P_{\gamma(t),p}\nabla_mk'(p) $$
\proofend

\subsection{Equivariant Kernel Functions}\label{sec:ekf}
Let
$$Q\colon\Omega^l(X)\to\Omega^k(X)$$
be an integral operator, i.e.
$$(Q\omega)_x=\int_X\hat{q}_{x,y}\omega_ydy$$
where $\hat{q}$ is a section in the exterior homomorphism bundle
$$\Hom(\Lambda^lTX^*,\Lambda^kTX^*)\to X\times X\ ,$$
i.e. $\hat{q}_{x,y}\in\Hom(\Lambda^lT_xX^*,\Lambda^kT_yX^*)\cong\Hom(\Lambda^kT_yX,\Lambda^l_xTX)$.
We can rewrite $\hat{q}$ in the form $\hat{q}\omega =\omega\circ q$ with some section
$$q\in\Gamma\Hom(\Lambda^kTX,\Lambda^lTX)$$
hence
$$(Q\omega)_x=\int_X \omega_y\circ q_{y,x}dy\ .$$
The operator $Q$ is $G$-equivariant if and only if
\[{eq:equivLL}q_{y,x}=g^{-1}q_{gy,gx}g\]
for all $x,y\in X$, $g\in G$.
Rank one symmetric spaces are two point homogeneous.
We can therefore map any pair $x,y\in X$ to the geodesic $\gamma$ by an isometry $g$ of $X$ so that
$$gy=p=\gamma(0)\qmbox{,}gx=\gamma(t)$$
where $t=d(x,y)\in\R^+_0$ is the distance of the two points.
The parallel transport satisfies the same equivariance \eqref{eq:equivLL} as $q$.
We can therefore rewrite $q=q^A$ in the form
\[{defLAt}q^A_{y,x}=g^{-1}q^A_{\gamma(0),\gamma(t)}g=g^{-1} A(t)P_{\gamma(0),\gamma(t)}g\]
with a function
\[{defAt}A\colon\R^+\to\Hom_H(\Lambda^kT_{\gamma(0)}X,\Lambda^lT_{\gamma(0)}X)=\Hom_H(\Lambda^k\m,\Lambda^l\m)\]
This establishes a one to one correspondence between equivariant kernel functions, respectively equivariant integral operators on one side and functions on $\R^+$ taking values in $\Hom_H(\Lambda^k\m,\Lambda^l\m)$ on the other,
{\small
\begin{align*}
\l(\Gamma\Hom(\Lambda^kTX,\Lambda^lTX)|_{X\times X\setminus\Delta X}\r)^G&\leftrightarrow C^\infty(\R^+,\Hom_H(\Lambda^k\m,\Lambda^l\m))\\
q^A&\leftrightarrow A
\end{align*}
}
In the case $l=k=0$ this is the bijection between radial kernel functions (see \cite{Szabo}) and functions on the positive real line.

\section{A Right Inverse for the Laplace Operator on Forms}\label{sec:rilof}
In this section we will construct an equivariant integral operator inverting the Laplacian.
Thus for $\omega\in\Omega^l(X)$ of compact support we look for a solution $\mu\in\Omega^l(X)$ of
$$\omega=\Delta\mu$$
in the form
\[{formofkernel}\mu_x=\int_X\omega_y\circ q^A_{y,x}dy\]
with an integral kernel determined by a function $A\in C^\infty(\R^+,\Hom_H(\Lambda^l\m,\Lambda^l\m))$ as in \eqref{defLAt}.
\begin{thm}\label{thm:main}
Let $\st{k_1,\ldots,k_n}\subset\k_0$ be the basis of \ref{ntt-basis}.
Let $h(t)$ be the mean curvature of the geodesic sphere in $X$ of radius $t$.
If $A\in C^\infty(\R^+,\Hom_H(\Lambda^l\m,\Lambda^l\m))$ is a solution of the ordinary second order differential equation
\[{maindeq}
\begin{split}A''(t)+h(t)A'(t)&=A(t)\RRR\\
&\hspace{-25mm}-\sum_{i=1}^n\frac{\ad_{k_i}^2A(t)-2\cosh(\lambda_i t)\ad_{k_i}A(t)\ad_{k_i}+\cosh(\lambda_i t)^2 A(t)\ad_{k_i}^2}{\sinh(\lambda_i t)^2/\lambda_i^2}
\end{split}\]
with initial condition
\[{eq:asympt}A(t)=\frac{t^{1-n}}{(n-1)\vol S^{n}}(1+o(1))\qmbox{as}t\to 0\ .\]
then $q^A$ as given in \eqref{defLAt} is a Green's function for the Laplacian.
Thus for any $\omega\in\Omega^l(X)$ with compact support we have
$$\Delta_x\int_X\omega_y\circ q^A_{y,x}dy=\omega_x\ .$$
\end{thm}
Note that the initial condition \eqref{eq:asympt} on $A$ is so that $\omega_p\circ q^A_{p,x}$ is asymptotic to the Newton Potential \eqref{eq:newton} in normal coordinates around $p$.

\proof We will show that if $A$ satisfies \eqref{maindeq} and \eqref{eq:asympt} then $\Delta_x\l(\omega_y\circ q^A_{y,x}\r)=\delta_y\omega_y$ for all $y$ and $\omega_y\in\Lambda^lT_yX^*$, i.e.
\[{eq:fsol}\int_X\scp{\Delta_x\alpha_x}{\omega_y\circ q^A_{y,x}}dx=\scp{\alpha_y}{\omega_y}\]
for all $\alpha\in\Omega^l(X)$ with compact support.
By homogeneity, it suffices to do this for $y=p$ and $x=\gamma(t)$ for some $t\in\R^+_0$.
To this end, let $\omega_p\in\Lambda^l\m^*=\Lambda^lT_pX^*$ and consider $\mu\in\Omega^l(X)$ given by
$$\mu_{e^{sk}\gamma(t)}(v_1\wedge\ldots\wedge v_l)=\omega_p\l(e^{sk}A(t)P_{p,\gamma(t)}e^{-sk}(v_1\wedge\ldots\wedge v_l)\r) $$
for $k\in\k$, $s\in\R$, $t\in\R^+_0$, $v_i\in T_{e^{sk}\gamma(t)}X$.
Thus $\mu$ is the composition
\[{eq:defomegaz}
\begin{split}
\mu_{e^{sk}\gamma(t)}&=\omega_p\circ q^A_{p,e^{sk}\gamma(t)}=\omega_p\circ e^{sk}\circ q^A_{p,\gamma(t)}e^{-sk}\\
&=\omega_p\circ e^{sk}A(t)P_{p,\gamma(t)}e^{-sk}\ .
\end{split}\]
We now compute $\Delta\mu$ at $\gamma(t)$.
From the Weitzenb\"ock formula \eqref{WBF},
\begin{align*}
\Delta\mu(\gamma(t))&=\l(-\nabla^2_{T,T}\mu-\sum_{i=1}^n\frac{1}{\norm{k_i}^2}\nabla^2_{k_i,k_i}\mu+\RR\mu\r)(\gamma(t))\\
&=-\omega_p\circ A''(t)P_{p,\gamma(t)}+\mu\circ\RRR(\gamma(t))\\
&\hspace{1cm}-\sum_{i=1}^n\frac{1}{\norm{k_i}^2}\nabla_{k_i}\nabla_{k_i}\mu(\gamma(t))+\sum_{i=1}^n\frac{1}{\norm{k_i}^2}\nabla_{\nabla_{k_i}k_i}\mu(\gamma(t))
\end{align*}
By \eqref{nablakk},
$$\sum_{i=1}^n\frac{1}{\norm{k_i}^2}\nabla_{k_i}k_i(\gamma(t))=-h(t)T(\gamma(t))\ ,$$
hence
\begin{align*}
\Delta\mu(\gamma(t))&=-\omega_p\circ (A''(t)+h(t)A'(t))\circ P_{p,\gamma(t)}\\
&\hspace{1cm}+\l(-\sum_{i=1}^n\frac{1}{\norm{k_i}^2}\nabla_{k_i}\nabla_{k_i}\mu+\mu\circ\RRR\r)(\gamma(t))
\end{align*}
In order to compute $\nabla_k\nabla_k\mu$ for $k=k_i\in\k_0$, we differentiate \eqref{eq:defomegaz} for $s$.
\begin{align*}
\DD{s}\l(e^{sk}A(t)\r)&=\dd{s}\l(e^{sk}A(t)\r)=\ad_kA\ .
\end{align*}
For the calculation of $\DD{s}P_{p,\gamma(t)}e^{-sk}$ let $u=\dd{r}u(r)\in T_{\gamma(t)}X$. Then
\begin{align*}
\l(\DD{s}P_{p,\gamma(t)}e^{-sk}\r)u&=\DD{s}\dd{r}P_{p,\gamma(t)}e^{-sk}u(r)\\
&=\DD{r}\dd{s}P_{p,\gamma(t)}e^{-sk}u(r)\\
&=-P_{p,\gamma(t)}\nabla_uk\\
&=-\cosh(\lambda t)\ad_k P_{p,\gamma(t)}u
\end{align*}
by \eqref{nablakinx}. Hence
$$\nabla_k\mu(\gamma(t))=\omega_p\circ\l(\ad_k A(t)-\cosh(\lambda t)A(t)\ad_k\r)\circ P_{p,\gamma(t)}\ .$$
Differentiating again in the same way, we get
$$\nabla_k\nabla_k \mu(\gamma(t))=\omega_p\circ Z_{t,k}^2(A(t))\circ P_{p,\gamma(t)}$$
where $Z_{t,k}(A(t))=\ad_k A(t)-\cosh(\lambda t)A(t)\ad_k$.
Thus
$$
\Delta\mu(\gamma(t))=-\omega_p\circ\l(A''(t)+h(t)A'(t)+\sum_{i=1}^n\frac{Z_{t,k_i}^2(A(t))}{\sinh(\lambda_i t)^2/\lambda_i^2}-A(t)\RRR\r)\circ P_{p,\gamma(t)}\ .\ 
$$
Therefore
$$\Delta_x\l(\omega_y\circ q^A_{y,x}\r)=0$$
holds for all $x=\gamma(t)\neq y=p=\gamma(0)$ and, by equivariance, for all $x,y$, $x\neq y$, provided that $A$ solves \eqref{maindeq}.
\par
For the proof of \eqref{eq:fsol} let $\alpha\in\Omega_c^l(X)$ be a $l$-form of compact support.
For $\rho>0$ let $\phi\in C^\infty(X)$ be a bump function such that
$$\phi(x)=\cases{1&x\in B_{\rho/2}(p)\\0&x\not\in B_{\rho}(p)}$$
and $|d\phi|<10/\rho$.
Then
\[{gfdsrtrtr}
\begin{split}\int_X\scp{\Delta_x\alpha_x}{\omega_p\circ q^A_{p,x}}dx&=\int_X\scp{\Delta_x(\phi\alpha)_x}{\omega_p\circ q^A_{p,x}}dx\\
&\hspace{5mm}+\int_X\scp{\Delta_x((1-\phi)\alpha)_x}{\omega_p\circ q^A_{p,x}}dx
\end{split}
\]
The right summand becomes
\begin{equation*}
\begin{split}\int_{X\setminus B_{\rho/2}(p)}&\scp{\Delta_x((1-\phi)\alpha)_x}{\omega_p\circ q^A_{p,x}}dx\\
&=\int_{X\setminus B_{\rho/2}(p) }\scp{((1-\phi)\alpha)_x}{\Delta_x\l(\omega_p\circ q^A_{p,x}\r)}dx=0
\end{split}
\end{equation*}
by what we proved above.
For the left summand, we compare the Laplace operator of $X$ with the euclidean Laplace operator $\Delta^0$ in normal coordinates near $p$ and then use that $\omega_p\circ q^A_{p,x}$ asymptotically coincides with $N$ which in turn is a fundamental solution to the euclidean Laplace operator.

To this end, let $e_i=\frac{\partial}{\partial x_i}$, $i=1\ldots d=n+1$, be the coordinate vector fields of normal coordinates near $p$.
The Riemannian connection $\nabla$ of $X$ near $p$ is then of the form
$\nabla_{v}=v+\Gamma_v$ with a local one form $\Gamma$ with values in the endeomorphism bundle of $TX$.
Since the coordinates are normal we have $\Gamma_v(p)=0$ for all $v\in T_pX$.
For the second covariant derivative we compute
$$\nabla^2_{v,v}=\nabla_v\nabla_w-\nabla_{\nabla_vw}=(D_v+\Gamma_v)(D_w+\Gamma_w)-D_{\Gamma_vw}-\Gamma_{\Gamma_vw}$$
$$=D^2_{v,w}+\Gamma_vD_w+\Gamma_wD_v-D_{\Gamma_vw}+\l(D_v\Gamma_w\r)+\Gamma_v\Gamma_w-\Gamma_{\Gamma_vw}$$
Thus for any $\alpha\in\Omega^l(X)$ we have
$$\nabla^2\alpha=D^2\alpha+F D\alpha +H\alpha$$
with tensor fields $F$ and $H$ and $F(0)=0$.
By the Weitzenb\"ock formula we also have such a formula for the Laplacians, i.e.
$$\Delta\alpha=\Delta^0\alpha+FD\alpha+H\alpha$$
with $F(p)=0$.
In the sequel for $x\in B_\rho(p)$, we will write $r=r(x)=d(x,p)$.
By our assumption we have
$$\omega_p\circ q^A_{p,x}=N\omega_p+R$$
with some $R\in\Omega^l(X)$ of order $R\sim r^{2-n}$.
\par
The left summand in \eqref{gfdsrtrtr} equals
\[{tyeala}=\int_{B_\rho(p)}\scp{\Delta_x(\phi\alpha)_x}{\omega_p\circ q^A_{p,x}}dx\ .\]
We will denote by $dx$ the volume form of the Riemannian metric from $X$ and by $dx^0$ the euclidean volume form on $B_\rho(p)$ induced by the normal coordinates.
If $g$ denotes the metric tensor expressed in the normal coordinates, then
$$dx=\sqrt{\det g}\ dx^0\sim dx^0$$
because $\det g(p)=1$.
The integral \eqref{tyeala} becomes
\begin{align*}
&=\underset{B_\rho(p)}{\int} \scp{\Delta^0_x(\phi\alpha)_x+F(x)D(\phi\alpha)_x+H(x)(\phi\alpha)_x}{N\omega_p+R}dx\\
&=\underset{B_\rho(p)}{\int}\scp{\underbrace{\Delta^0_x(\phi\alpha)_x}_{\rho^{-2}}+\underbrace{F(x)D(\phi\alpha)_x}_{\rho\times\frac{1}{\rho}=\rho^0}+\underbrace{H(x)(\phi\alpha)_x}_{\rho^0}}{\underbrace{N\omega_p}_{\rho^{1-n}}+\underbrace{R}_{\rho^{2-n}}}\underbrace{\sqrt{\det g}}_{1}\ dx^0
\end{align*}
In the limit $\rho\to0$ the terms involving $F$, $H$ and $R$ vanish and thus for the integral \eqref{tyeala} we finally get
\begin{align*}
&=\int_X\scp{\Delta_x\alpha_x}{\omega_p\circ q^A_{p,x}}dx\\
&=\lim_{\rho\to 0}\int_{B_\rho(p)}\scp{\Delta^0_x(\phi\alpha)_x}{N\omega_p}dx^0\\
&=\lim_{\rho\to 0}\scp{\phi(p)\alpha_p}{\omega_p}=\scp{\alpha_p}{\omega_p}
\end{align*}
since $N$ is a fundamental solution for the euclidean Laplacian.
\proofend
\subsection{Existence of a solution}
A solution $A$ of the equation \eqref{maindeq} in Theorem \ref{thm:main} always exists.
To see this we rewrite the differential equation \eqref{maindeq} in the form
$$(\sigma A')'= W(\sigma A)$$
where $WA$ denotes the right hand side of \eqref{maindeq}.
Substituting
$$z=\frac{\sigma A}{t}\qmbox{and}y=\sigma A'$$
this becomes
\begin{align*}
z'&=-\frac{z}{t}+\frac{\sigma A'}{t}+\frac{\sigma' A}{t}=\frac{1}{t}\l(y+(th-1)z\r)\\
y'&=W\sigma A=tW z
\end{align*}
From \eqref{eq:sigmah} and \eqref{maindeq} we can write $h$ and $W$ in the form
$$h(t)=\frac{n}{t}+th_0(t)\qmbox{and}W(t)=\frac{W_{-2}}{t^2}+W_0(t)$$
with even analytic functions $h_0$ and $W_0$.
We thus have the regular singular first order differential equation
\begin{align*}
\pmx{z\\y}'&=\frac{1}{t}\pmx{th(t)-1&1\\t^2W(t)&0}\pmx{z\\y}\\
&=\l(\frac{1}{t}\pmx{n-1&1\\W_{-2}&0}+t\pmx{h_0(t)&0\\W_0(t)&0}\r)\pmx{z\\y}\ .
\end{align*}
We now substitude $t=e^{-s}$, $Y(s)=(z,y)(e^{-s})$, $Y'=-e^{-s}(z,y)'(e^{-s})$ to transform this into
\[{eq:levinson}Y'(s)=\l(-\pmx{n-1&1\\W_{-2}&0}-e^{-2s}\pmx{h_0(e^{-s})&0\\W_0(e^{-s})&0}\r)Y(s)\ .\]
Since $W_{-2}$ is Hermitian and nonnegative, we may assume that $W_{-2}$ is a diagonal matrix with nonnegative eigenvalues $w_i$, $i=1,\ldots,\dim\Hom_H(\Lambda^l\m,\Lambda^l\m)$.
We can therefore diagonalise the matrix $\pmx{n-1&1\\W_{-2}&0}$ with $(2\times 2)$-blocks of the form $\pmx{n-1&1\\w_i&0}$.
Since the $w_i$ are nonnegative, these blocks have 2 different eigenvalues.
It follows that they all can be diagonalised.

Therefore the equation \eqref{eq:levinson} satisfies the conditions of the Levinson Theorem for asymptotic constant (diagonalisable) coefficient systems.
\begin{thm}[Theorem 1.8.1 in \cite{Eastham}]\label{levinson}
Let $C$ be a diagonalizable $m\times m$ matrix and $R(s)$ a matrix valued function, $s\in\R$.
Assume that
\[{levinson-cond}\int_a^\infty\norm{R(s)}ds<\infty\]
for all $a\in\R$.
If $\set{u_k}{1\leq k\leq m}$ are linearly independent eigenvectors of $C$ corresponding to the eigenvalues $\lambda_k$, $1\leq k\leq m$, then the differential equation
\[{levinson-deq}Y'(s)=\l(C+R(s)\r)Y(s)\]
has solutions with the asymptotic form
$$Y_k(s)=\l(u_k+o(1)\r)e^{\lambda_ks}\qmbox{as}s\to\infty\ .$$
\end{thm}
In particular, for the eigenvalue $0$ of $C=\pmx{n-1&1\\-W_{-2}&0}$, we get a solution $Y(s)=\pmx{\frac{-1}{n-1}\id\\\id}+o(1)$ for \eqref{eq:levinson} which corresponds to a solution $A$ for \eqref{maindeq} with
$$\sigma(t)A'(t)=\id+o(1)\qmbox{,}\frac{\sigma(t)A(t)}{t}=\frac{-\id}{n-1}+o(1)\ .$$
Hence, because $\sigma(t)=t^n(\vol S^n+o(1))$,
$$A(t)=\frac{-t}{(n-1)\sigma(t)}(\id+o(1))=\frac{-t^{1-n}}{(n-1)\vol{S^n}}(\id+o(1))$$
as $t\to 0$.

The operators $d$, $d^*$ and $\Delta$ preserve this class of kernels.
For later reference, we note their action on the corresponding radial function $A$.

\begin{lm} Let $A\in C^\infty(\R^+,\Hom_H(\Lambda^k\m,\Lambda^l\m))$ and $Q^A$ be the integral operator with kernel $q^A$, i.e.
$$(Q^A\omega)_x=\int_X \omega_y\circ q^A_{y,x}dy$$
for $\omega\in\Omega^l(X)$.
Then $dQ^A=Q^{dA}$ and $d^*Q^A=Q^{\dsa}$ with
$$dA\in C^\infty(\R^+,\Hom_H(\Lambda^{k+1}\m,\Lambda^l\m))\ ,$$
$$\dsa\in C^\infty(\R^+,\Hom_H(\Lambda^{k-1}\m,\Lambda^l\m))\ ,$$
given by
$$dA(t)=T\wedge A'(t)+\sum_{i}\frac{m_i(p)\wedge\l(\ad_{k_i}A(t)-\cosh(\lambda_i t)A(t)\ad_{k_i}\r)}{\sinh(t\lambda_i)/\lambda_i}\ ,$$
$$-\dsa(t)=T\neg A'(t)+\sum_{i}\frac{m_i(p)\neg\l(\ad_{k_i}A(t)-\cosh(\lambda_i t)A(t)\ad_{k_i}\r))}{\sinh(t\lambda_i)/\lambda_i}\ .$$
\end{lm}
\proof This follows from \eqref{diffcodiff} by the same calculation as in the proof of Theorem \ref{thm:main}.
\proofend

\subsection{Hyperbolic Space}\label{hypsp}
In space forms $X^m$ the orthogonal complement $\m_0$ of the unit vector $T\in T_pX=\m$ is irreducible for $H$. We have a decomposition 
$$\Lambda^l\m=T\wedge\Lambda^{l-1}\m_0\oplus\Lambda^l\m_0 $$
as modules for $H$. For $A\in C^\infty(\R^+,\Hom_H(\Lambda^l\m,\Lambda^l\m))$ we can therefore make the ansatz
$$A(t)(T\wedge\mu_0+\mu_1)=\alpha(t)T\wedge\mu_0+\beta(t)\mu_1$$
with
$$\mu_0\in\Lambda^{l-1}\m_0\ ,\ \mu_1\in\Lambda^l\m_0\qmbox{and}\alpha,\beta\in C^\infty(\R^+,\R)\ .$$
Let $k=m-l=n+1-l$. A straightforward calculation shows that in hyperbolic $m$-space the curvature endomorphism $\RRR$ on $\Lambda^l\m$ from the Weitzenb\"ock formula is multiplication with $-kl$ and
$$\ad_{k_i}m_i=-T\qmbox{,}\ad_{k_i}T=m_i\qmbox{and}\ad_{k_i}m_j=0\qmbox{if}j\neq i\ .$$
The differential equation \eqref{maindeq} for $A(t)$ becomes

\[{hypspacer4}
\begin{split}
\alpha''(t)+n\coth(t)\alpha'(t)&=\\
&\hspace{-20mm}-kl\alpha(t)+k\frac{1+\cosh(t)^2}{\sinh(t)^2}\alpha(t)-2k\frac{\cosh(t)}{\sinh(t)^2}\beta(t)
\\
\beta''(t)+n\coth(t)\beta'(t)&=\\
&\hspace{-20mm}-kl\beta(t)+l\frac{1+\cosh(t)^2}{\sinh(t)^2}\beta(t)-2l\frac{\cosh(t)}{\sinh(t)^2}\alpha(t)
\end{split}
\]
The solution to \eqref{hypspacer4} satisfying the initial condition \eqref{eq:asympt} will be Laurent series $\alpha,\beta$ starting with $\frac{t^{1-n}}{(n-1)\vol S^{n}}$, whose higher order terms can be computed from \eqref{hypspacer4}.
The expressions of these solutions to \eqref{hypspacer4} in terms of elementary functions are quite complicated. For example, in the case of hyperbolic $3$-space and the Laplacian on $1$-forms, $l=1$, $k=2$, we get solutions
\begin{align*}
\alpha(t)&=\frac{(2t+1)\sinh(t)-(t^2+t)\cosh(t)}{4\pi\sinh(t)^3}\ ,\\
\beta(t)&=\frac{3\sinh(t)^2-(\cosh(t)+2t)\sinh(t)+t^2+t}{8\pi\sinh(t)^3}-\frac{t}{4\pi(\cosh(t)+1)}
\end{align*}
which have also been found in \cite{dtg}.
In fact, $\alpha,\beta$ above are given by $\alpha(t)=-\phi_2(t)-\phi_3''(t)$ and $\beta(t)=-\phi_2(t)-\phi_3'(t)/\sinh(t)$
where the functions $\phi_2$, $\phi_3$ are as in \cite{dtg}, Theorem 3, (3) page 20.


\section{The Biot-Savart Formula and The Linking Number}\label{sec:biotsavart}
By Hodge Theory, on any oriented compact Riemannian manifold, a right inverse $Q$ of the Laplace operator on differential forms provides a right inverse to the Cartan differential $d$ and the codifferential $d^*$ simultaneously.
Thus let $Q$ be a right inverse of the Laplacian, i.e. $\Delta Q \omega=\omega$ for all $l$-forms $\omega$ in the image of $\Delta$.
If, in particular, $\alpha\in\Omega^l(X)$ is exact, then $\alpha=dd^*\!Q\alpha$ and $d^*\!d^*\!Q\alpha=0$.
Following \cite{ctg}, \cite{dtg} and \cite{tk} we call an operator $B$ so that $dB\alpha=\alpha$ and $d^*\!B\alpha=0$ a Biot-Savart operator.
This is motivated by the analogy of these equations with the Maxwell equations for the magnetic field $B$ caused by a stationary current $\alpha$.

The differential equation \eqref{maindeq} is invariant under complex rescaling of the metric on $X$.
If $A(t)$ solves \eqref{maindeq} for the symmetric space whose Jacobi operator has eigenvalues $\lambda_1^2,\ldots,\lambda_n^2$, then for $\alpha\in\C\setminus\st{0}$,
$$A_\alpha(t)=\alpha^{n-1}A(\alpha t)$$
solves \eqref{maindeq}/\eqref{eq:asympt} with $\lambda_i$ replaced by $\alpha\lambda_i$ and the curvature tensor $R$ replaced by $\alpha^2 R$.
In particular, $A_i(t)=i^{n-1}A(it)$ solves \eqref{maindeq} for the compact dual $X_c$ of $X$.
Since $A$ is holomorphic on $\C\setminus\st{0}$, $A_\alpha(t)$ is holomorphic in $\alpha$.

On the compact Riemannian manifold $X_c$ Hodge theory gives an $L^2$-orthogonal splitting of $\Omega^l(X_c)$ as in the diagram,
\begin{center}
\begin{tikzcd}[column sep=1, row sep=15]
\Omega^l(X_c)\arrow[bend left=20]{dd}{d} & = & H^l(X_c)\arrow{d}{d} & \oplus & d\Omega^{l-1}(X_c)\arrow{dl}{d} & \oplus & d^*\Omega^{l+1}(X_c)\arrow[bend left=20]{lldd}[near start]{d}[near end, description]{\cong}\\
&&0&0&&&&0\\
\Omega^{l+1}(X_c)\arrow[bend left=20]{uu}{d^*} & = & H^{l+1}(X_c)\arrow{u}{d^*} & \oplus & d\Omega^{l}(X_c)\arrow[bend left=20]{rruu}[near start]{d^*}[near end, description]{\cong} & \oplus & d^*\Omega^{l+2}(X_c)\arrow{ur}{d^*}
\end{tikzcd}
\end{center}
The Cartan differential $d$ and its adjoint restrict to isomorphisms as indicated.
The Laplacian $\Delta_c=dd^*+d^*d$ of $X_c$ preserves the above splitting, vanishes on $H^l(X_c)=\ker\Delta_c=\ker d\cap\ker d^*$ and is isomorphic on both $d\Omega^{l-1}$ and $d^*\Omega^{l+1}$.

Let $\Delta_c^{-1}$ be the extension of the inverse of this isomorphism by $0$ to all of $\Omega^l(X_c)$.
By equivariance, the integral kernel of $\Delta_c^{-1}$ must be of the form \eqref{formofkernel} with a function $\tilde{A}_c$ satisfying \eqref{maindeq} and \eqref{eq:asympt} for the compact rank one symmetric space $X_c$.
We can therefore choose a solution $A$ to \eqref{maindeq} with $A_i=\tilde{A}_c$.
In this section, $A$ will always be such a solution of \eqref{maindeq}/\eqref{eq:asympt}.

Let $\omega\in d\Omega^l_c(X)$.
Then $\omega=dd^*\Delta_c^{-1}\omega=dQ^{\dsa_i}\omega=dQ^{\dsa_{i\alpha}}\omega$ for all $\alpha\in\R^+$.
Since this is holomorphic in $\alpha$, the identity holds for all $\alpha\in\C\setminus\st{0}$.
Thus $\omega=dQ^{\dsa}\omega$, for the above choice of $A$.
We thus have the following formula for the Biot-Savart operator on non-compact rank one symmetric spaces.
\begin{thm}\label{thm:bs}
Let $X$ be a noncompact rank one symmetric space and let $\omega\in\Omega_c^l(X)$ be a closed $l$-form of compact support.
Then
$$B=d^*\Delta^{-1}\omega=\int_X\omega_y\circ q^{\dsa}_{y,x}dy$$
satisfies
\[{biotsavartproblem}dB=\omega\qmbox{and}d^* B=0\ .\]
\end{thm}

There also is an integral formula for the linking number similar to the classical Gauss formula for the linking number of loops in $\R^3$.
\begin{thm}\label{thm:gs} Let $X^d$ be a noncompact rank one symmetric space.
Let $K^k,L^l\subset X^d$, $k+l+1=d$ be disjoint closed submanifolds of dimensions $k$ and $l$ respectively.
Let $A$ be the kernel function inverting the Laplace Beltrami operator on $k$-forms on $X$.
Then the linking number of $K$ and $L$ is given by
$$L(K,L)=\int_K\int_L\vol^L(y)\wedge q^{\dsa}_{y,x} \vol^K(x)dydx$$
where $\vol^L(y)$ and $\vol^K(x)$ denote the metric volume elements of $L$ at $y\in L$ and of $K$ at $x\in K$ respecively.
\end{thm}
\proof
Let $K^k,L^l\subset X^d$, $k+l+1=d$ be disjoint closed submanifolds of dimensions $k$ and $l$ respectively in a $d$-dimensional manifold $X$.
Let $i\colon K\hookrightarrow X$ and $j\colon L\hookrightarrow X$ be the inclusions.
We assume that $H_c^{d-k}(X)\cong H_k(X)=0=H^{d-l}_c(X)\cong H_l(X)$.
Following the notation of \cite{botu} we let $\eta_K\in\Omega^{d-k}_c(X)$ be the Poincar\'e dual of $K$, i.e. a compactly supported form on $X$ such that
$$d\eta_K=0\qmbox{and}\int_K i^*\alpha=\int_X\alpha\wedge\eta_K$$
for all $\alpha\in\Omega^k(X)$ with $d\alpha=0$.
Then the linking number of $K$ and $L$ in $X$ is
$$L(K,L):= \int_X\beta\wedge\eta_K=\int_Ki^*d^{-1}\eta_L$$
where $\beta$ is any form of degree $d-l=k+1$ with $d\beta=\eta_L$, for instance $\beta=d^*\Delta^{-1}\eta_L=Q^{\dsa}\eta_L$,
$$\beta_x=\int_X\eta_{L,_y}\circ q^{\dsa}_{y,x}dy\ ,$$
we obtain
$$L(K,L)=\int_Ki^*d^*\Delta^{-1}\eta_L=\int_K\int_X\eta_{L,y}\l(q^{d^* A}_{y,x}\vol^K(x)\r) dy dx$$
where $\vol^K(x)$ denotes the metric volume element of $K$ at $x$.
Let $NL$ be some tubular neighbourhood of $L$ induced by normal coordinates.
Thus we choose $r>0$ so small that the Riemannian exponential map
$$D_r\nu(L,X)\to NL\qmbox{,}v\in\nu_y(L,X)\mapsto\exp_yv$$
is a diffeomorphism on the $r$-disk bundle of the normal bundle $\nu(L,X)$ of $L$ in $X$.
For $1\geq\rho>0$, let $N_\rho L$ denote the tubular neighbourhood of $L$ obtained from $NL$ be rescaling with $\rho$, that is
$$N_\rho L=\set{\rho v}{v\in NL}\subset NL\ .$$
To represent the Poincar\'e dual of the submanifold $L$ we may take the extension by $0$ of the Thom form of the normal bundle of $L$.
This is a closed $(k+1)$-form on $NL$ vanishing near the bondary of $NL$ and such that its integral over the fibres of $NL$ is always $1$.
Let $\chi^\rho$ be the diffeomorphism
$$\chi_\rho\colon N_\rho L\to NL\qmbox{,}v\mapsto\frac{1}{\rho}v$$
and $\eta_L^\rho\in\Omega(X)$ be the extension by zero of 
$$\eta_L^\rho:=\chi_\rho^*\eta_L\ .$$
The forms $\eta_L^\rho$ are all closed.
Since the diffeomorphisms $\chi_\rho$ preserve the fibres of the tubular neighbourhoods, the integral over the fibres of $NL\to L$ of $\eta_L^\rho$ is always $1$.
Thus the forms $\eta_L^\rho$ are all Thom forms for the normal bundle, or, Poincar\'e duals of $L$.
For the linking number it follows that
\[{eq:kf57}L(K,L)=\int_K\int_{N_\rho L}\eta^\rho_{L,y}\l(q^{\dsa}_{y,x}\vol^K(x)\r) dy dx\]
for all $1\geq\rho>0$.
We now take the limit $\rho\to 0$.
Using the Riemannian metric, we split the tangent bundle of $NL$ into vertical and horizontal distribution.
The Thom form $\eta_L$ can then be written as a sum
$$\eta_L=\omega_0+\sum_{i=1}^l\alpha_i\wedge\omega_i$$
with vertical forms $\omega_i$ of degree $k+1-i=d-l-i$ and horizontal forms $\alpha_i$ of degree $i$.
$\omega_0$ is a volume form of the fibres.
As $\rho\to 0$,
$$(\chi^\rho)^*\l(\sum_{i=1}^l\alpha_i\wedge\omega_i\r)\to 0\qmbox{in}L^1$$
because the differential of $\chi^\rho$ on $\omega_i$ scales with a factor of $1/\rho^{k+1-i}$ but the support scales with $\rho^k$.
Taking the limit $\rho\to 0$ in \eqref{eq:kf57} gives
\[{eq:kf533}L(K,L)=\lim_{\rho\to 0}\int_K\int_{N^\rho L}(\chi^\rho)^*\omega_{0,y}\l(q^{\dsa}_{y,x}\vol^K(x)\r)\ dydx\]
$$=\int_K\int_L\widetilde{\omega}_{0,y}\l(q^{\dsa}_{y,x}\vol^K(x)\r) dydx$$
where $\widetilde{\omega}_{0}\in\Gamma\l(\Lambda^{k+1}TX|_L^*\r)$ is such that $\widetilde{\omega}_{0,y}(v_0,v_1,\ldots v_k)=1$ if $(v_0,v_1,\ldots v_k)$ is an oriented orthornormal basis for the orthogonal complement of $T_yL$ in $T_yX$ and $0$ if one of the $v_i$ is tangent to $L$.
Thus the above integral becomes
$$=\int_K\int_L\vol^L(y)\wedge q^{\dsa}_{y,x}\vol^K(x)\ dydx\ .$$
\proofend


\end{document}